\documentclass[article]{amsart}
\usepackage{cases}
\usepackage{amsmath, amsfonts, pifont, amssymb}
\usepackage{verbatim}
\usepackage[bookmarks=true]{hyperref}
\usepackage{geometry}
\usepackage{color}
\newtheorem{theorem}{Theorem}[section]
\newtheorem{lemma}[theorem]{Lemma}

\newtheorem{assume}[theorem]{Assumption}

\newcommand{\beq}{\begin{equation}}
\newcommand{\eeq}{\end{equation}}
\newcommand{\beqq}{\begin{equation*}}
\newcommand{\eeqq}{\end{equation*}}

\newcommand{\PPP}{\mathcal{P}}
\theoremstyle{definition}

\theoremstyle{remark}
\newtheorem{remark}[theorem]{Remark}

\numberwithin{equation}{section}

\numberwithin{equation}{section}

\begin{document}
\title{A note on minimal mass blow up for inhomogeneous NLS}

\author{ Chenjie Fan}
\address{Academy of Mathematics and Systems Science and Hua Loo-Keng Key Laboratory of Mathematics,}
\address{Chinese Academy of Sciences, Beijing, China.}
\email{fancj@amss.ac.cn}

\author{Shumao Wang}
\address{Academy of Mathematics and Systems Science, The Chinese Academy of Sciences, Beijing 100190, CHINA.}
\email{wangshumao@amss.ac.cn}
\maketitle
\begin{abstract}
We revisit the work of \cite{raphael2011existence} on the minimal mass blow up solution for $iu_{t}+\Delta u=-k(x)|u|^{2}u$, and extend the construction of such a solution to the $k\in C^{2}$ case. 
\end{abstract}

\section{Introduction}
\subsection{Statement of main results and motivations}
Consider inhomogeneous NLS in $\mathbb{R}^{2}$,
\begin{equation}\label{eq: knls}
iu_{t}+\Delta u=-k(x)|u|^{2}u,\\
u(0,x)=u_{0}\in H^{1}(\mathbb{R}^{2}),
\end{equation}

We assume 
\begin{assume}\label{as: main}
	$k$ is $C^{2}$ with
\begin{equation}
k(0)=1, 0\leq k(x)\leq 1, det\nabla^{2}k(0)\neq 0,
\end{equation}
and we assume that $k$, derivative of $k$, and second derivative of $k$ satisfy 
\begin{equation}
\|k\|_{L_{x}^{\infty}}+\|\nabla k\|_{L_{x}^{\infty}}+\sum_{i,j=1,2}\|{\partial_{ij}k}\|_{L_{x}^{\infty}}<\infty
\end{equation}
\end{assume}

It can be proven for all initial data with $\|u_{0}\|_{L_{x}^2}<\|Q\|_{L_{x}^{2}}$, the associated solution is global. Here, the ground state $Q$ is the unique positive radial solution which solves 
\begin{equation}\label{eq: groundstate}
-\Delta Q+Q=|Q|^{2}Q,
\end{equation}

In the current article, we construct finite time blow up solution to \eqref{eq: knls} with mass equal to $\|Q\|_{L_{x}^{2}}$. Such solutions are sometimes called minimal mass blow up solutions. 
\begin{theorem}\label{thm: main}
Under Assumption \ref{as: main}, there exists initial data $u_{0}\in H^{1}$ with $\|u_{0}\|_{L_{x}^{2}}=\|Q\|_{L_{x}^{2}}$, such that the associated solution blows up in finite time.
\end{theorem}

\begin{remark}
It is expected one can get rid of the assumption $det\nabla^{2}k(0)\neq 0$, but we didn't handle it in the current article.
\end{remark}

Since (constructive) blow up is local in space in some sense, we further assume WLOG that $k(x)<1$, $x\neq 0$, for convenience, thus all possible minimal blow up solution must concentration at the origin.

This work has been motivated by the seminal work \cite{raphael2011existence}, where minimal mass blow up solutions to \eqref{eq: knls} has been constructed for $k\in C^{5}$. The construction of for $k$ with lower regularity has also been studied. In \cite{banica2010minimal}, the construction has been done under the assumption that $k\in C^{4}$, with degenerate conditions $\partial_{ij}k(0)=0, i,j=1,2$. Minimal mass of blow up solution has also been constructed in \cite{matsui2021minimal}, for $k\in C^{2}$ but with a degenerate condition $|k(x)-1|\lesssim |x|^{2+\delta},\delta>0$ for $x$ close to zero. Morally speaking, it is like assuming $k\in C^{2+\delta}$ with $\partial_{ij}k(0)=0, i,j=1,2$.

Note the $C^{2}$ regularity is in some sense critical in this problem. Indeed, it has been proved in \cite{merle1996nonexistence}, that if one has 
\begin{equation}\label{eq: condi1}
x\nabla k(x)\leq -|x|^{1+\alpha}, \alpha\in (0,1), \text{ for } |x|\ll 1.
\end{equation}
then no minimal mass blow up solution can exist for \eqref{eq: knls}. Note that \eqref{eq: condi1} implies $k$ cannot be $C^{2}$.

It has also been observed and  noted in \cite{raphael2011existence}, that the impact caused by $\nabla^{2}k(0)$ is a priori non-perturbative \footnote{But it turns out still be perturbative in the sense there is an extra algebraic cancellation which makes this part does not impact the blow up rate of the minimal mass blow up solution}. Thus, to study $k$ at exact $C^{2}$ regularity is in some sense subtle. Also, from this perspective, it may be expected Theorem \ref{thm: main} can be obtained without the $det \nabla^{2}k(0)\neq 0$, condition.

One will see the whole proof scheme of current article follows \cite{raphael2011existence} closely, but a modified Lyapounov function will be needed to over come the low regularity issue for $k$.

It should be noted that in the work of \cite{raphael2011existence}, the uniqueness of the minimal mass blow up solutions has been obtained, this part seems to depend on the regularity of $k$ in a more sensitive way. We didn't investigate this issue in the current article.

We also take this chance to mention one can apply the method in \cite{raphael2011existence} to construct Bourgain-Wang solutions on $\mathbb{T}^{2}$.
\subsection{Background for NLS and constructive blow up}
Focusing (2D) Mass critical NLS is one of the most typical Model in nonlinear Dispersive PDEs, in particular the study of blow up problems.
\begin{equation}\label{eq: nls}
\begin{cases}
iu_{t}+\Delta u=-|u|^{2}u,\\
u(0,x)=u_{0}(x)\in H^{1}(\mathbb{R}^{2})
\end{cases}
\end{equation}
It has three conservation laws\footnote{We leave the $M$, $E$ notation for the energy and mass for \eqref{eq: knls}.},
\begin{enumerate}
	\item Mass: ${M}(u_{0})=\int |u_{0}|^{2}$.
	\item Energy: ${E}(u_{0})=\frac{1}{2}\int |\nabla u_{0}|^{2}-\frac{1}{4}\int |u_{0}|^{4} $
	\item Momentum: ${P}(u_{0})=\Im \int u_{0}\nabla \bar{u}_{0}$
\end{enumerate}	
Note that for focusing NLS, the energy can be negative. 

Its symmetry includes translation in space, and time, scaling (which is $L^{2}$ invariant), and phase, and pseudo-conformal symmetry. The last one states that if $u$ solves \eqref{eq: nls}, so does
\begin{equation}\label{eq: pesudo}
\frac{1}{T-t}u(\frac{1}{T-t},\frac{x}{T-t})e^{-i\frac{|x|^{2}}{4(T-t)}}
\end{equation}
where $T>0$ is a parameter.

Observing $Q(x)e^{it}$ solves \eqref{eq: nls} if $Q$ solves \eqref{eq: groundstate}, and applying Pesudo-conformal symmetry to $Q(x)e^{it}$ gives explicit minimal blow up solutions,
\begin{equation}
S(t,x):=	\frac{1}{T-t}Q(\frac{x}{T-t})e^{-i\frac{|x|^{2}}{4(T-t)}}e^{i\frac{1}{T-t}}.
\end{equation}

 It has been proved, that for all $L^{2}$ solutions, with mass below ground state, are global and scattering, \cite{weinstein1983nonlinear}, \cite{dodson2011global}. It has been further prove that minimal mass blow up solution to \eqref{eq: nls} are  $S(t,x)$ up to natural symmetry, \cite{merle1993determination},\cite{dodson2023determination}, \cite{dodson2021determination}.

Pesudo conformal symmetry is in some sense specific to mass critical model, it is related to the Virial identify, \cite{glassey1977blowing}, which gives existence of many blow up solutions to \eqref{eq: nls},
\begin{equation}\label{eq: virial}
\partial_{tt}\int|x|^{2}|u|^{2}=16{E}(u).
\end{equation}

In some sense, both pseudo-conformal symmetry and Virial identity requires extra integrability of initial data, i.e. $xu_{0}\in L^{2}$, and a lot effort has been made historically to get rid of this integrability, see in particular \cite{merle1993determination},\cite{merle2005blow}. There are also a lot of efforts to find a more dynamic alternative for the pesudo-conformal symmetry. And this is indeed one of the main motivation of \cite{raphael2011existence} to study model \eqref{eq: knls}. It should be noted that \eqref{eq: knls} is structurally similar to \eqref{eq: nls}. It also admits mass and energy conservation\footnote{Note that the energy of \eqref{eq: nls} and \eqref{eq: knls} are slightly different.}, with
\begin{enumerate}
	\item $M(u_{0}):=\int |u_{0}|^{2},$
	\item $E_{in}(u_{0}):=\frac{1}{2}\int |\nabla u_{0}|^{2}-\frac{1}{4}\int k(x)|u_{0}|^{4} $.
\end{enumerate}

But the in-homogeneous model \eqref{eq: knls} does not admit explicit pseudo-conformal symmetry. We note also does not have exact scaling and (space) translation invariance, also not conservation of momentum.

Another typical model to study NLS without exact pseudo-conformal symmetry is NLS on torus, in particular . For example, in \cite{ogawa1990blow}, one applies the strategy in \cite{merle1990construction} to construct minimal mass blow up solutions in the periodic case. We note that a prioi it is not even clear there exists blow up, since $x$ is not a well defined function on torus, \eqref{eq: virial} is not directly applicable. 

The current work is in the field of constructive blow up. It aims to understand the blow up mechanism of solutions from a constructive view point. For comparison, Virial type argument gives the existence of many blow up solutions but does not say too much how the solution will actually blow up.

  It is a very active and vast research field, It is impossible to do an exhaustive review here, and we briefly review it focusing on NLS related problems. It is closely related to the modulation analysis, which aims to reduces the study of PDEs to finite dimensional ODEs, \cite{weinstein1985modulational}. One very influential work in the field is Merle and Rapha\"el's series work on log-log blow up solutions to mass critical NLS, \cite{merle2003sharp},\cite{merle2005blow},\cite{merle2006sharp},\cite{merle2004universality}, they highlight, among others,
  \begin{enumerate}
  \item One may center the solution not around exact ground state $Q$, but a modified version of $Q_{b}$, (where $b$ is a parameter). 
  \item It is not enough to perform energy estimates (based on energy conservation type computations), but also important to include the information from Virial identity
  \item Based on the second point, the choice of orthogonality condition ( and study of certain coerciveness) should not be purely based on the linear perturbation around ground state $Q$, but also based on certain linearized operator related to Virial.
  \end{enumerate}

See also \cite{perelman2001blow} for the first mathematical construction of log-log blow up. 

Another important idea goes back to \cite{merle1990construction}, which is one of most early constructive blow up for NLS. Roughly speaking, to build a blow up solution $u$ to \eqref{eq: nls} with a asymptotic structure $S$ which blows up at time $T$, rather than perturb $S(T-\epsilon)$ for $\epsilon\ll 1$ to make it blow up like $S$, one may also study \eqref{eq: nls} backwards for $u_{\epsilon}(T-\epsilon)=S(T-\epsilon)$,( and let $u_{\epsilon}$ be the associated solution), and prove that $u_{\epsilon}(T-t_{0})$ admits some uniform bound, thus one can use compactness argument to find a limit of $u_{\epsilon_{n}}(T-t_{0})$, which give rise the desired blow up solution.

In some sense, the seminal work \cite{raphael2011existence}, (which inspires our work), combines the above two ideas and extends it naturally in several way, in particular, if one can construct a modified ground $Q_{b}$ for a single parameter $b$, one can also construct a modified ground state $Q_{\PPP}$ for a family of parameters $\PPP$.

We refer to \cite{perelman2001blow},\cite{merle1990construction},\cite{merle2003sharp},\cite{merle2005blow},\cite{merle2006sharp},\cite{planchon2007existence},\cite{merle2004universality},
\cite{fan2017log},\cite{fan2023construction},\cite{rockner2024multi},\cite{bourgain1997construction},\cite{holmer2012blow}, \cite{rockner2024multi},\cite{bourgain1997construction},\cite{burq2003two}, \cite{raphael2006existence}, \cite{raphael2011existence},\cite{ogawa1990blow} and reference therein for constructive blow up related to log-log and Bourgain-Wang (or $S(t)$ like ) blow up for mass critical NLS.
 We also want to mention \cite{martel2018strongly} which gives new blow up rate for mass critical NLS by considering strongly interacting multi-bubbles.
 
We finally mention recent breakthrough for defocusing energy supercritical NLS and compressible fluid blow up solutions by Merle, Rapha\"{e}l, Rodnianski and 	Szeftel in \cite{merle2022implosion1}, \cite{merle2022blow} and \cite{merle2022implosion}, and see also the work \cite{buckmaster2022smooth} and \cite{cao2023non}.

\subsection{Structure of the article}
We present preliminary for modulation theory in Section \ref{sec: pre}. Then we set up the main bootstrap in Section \ref{sec: mainans}, and we prove the bootstrap Lemma in Section \ref{sec: bootstrap}.
\subsection{Notation}
We let 

\[
<f,g>:=\int_{\mathbb{R}^{2}}f(x)\overline{g(x)}dx.
\]

And $A\lesssim B$ means, there exists a universial constant $C$
such that 
\[
A\leq CB.
\]

\section{Acknowledgement}
We thank professor Ping Zhang for discussion and encouragement. This work has been partially supported by National Key R\&D Program of China, 2021YFA1000800, CAS Project for Young Scientists in Basic Research, Grant No. YSBR-031, and NSFC
Grant No. 12288201.
\section{Preliminary}\label{sec: pre}
\subsection{The spectral structure}
\subsubsection{Basics}
We recall  Schr\"{o}dinger operator $L_{+}, L_{-}$, which arise naturally when one linearized around ground state $Q$.
\begin{equation}
L_{+}:=-\Delta+1-3Q^{2}, \quad L_{-}:=-\Delta+1-Q^{2}.
\end{equation}
Formally speaking, let $J(f):=-\Delta f+f+|f|^{2}f$, (thus $J(Q)=0$), one has 
\begin{equation}
J(Q+\epsilon)=L_{+}\epsilon_{1}+iL_{-}\epsilon_{2}+O(\epsilon^{2})	
\end{equation}

The spectral property of $L_{\pm}$ has been well studied in \cite{weinstein1985modulational}, and one has 
\begin{equation}\label{eq: spectralba1}
L_{+}(\nabla Q)=0, \quad L_{+}\Lambda Q=-2Q, \quad L_{+}\rho=|y|^{2}Q,
\end{equation}
\begin{equation}\label{eq: spectralba2}
L_{-}Q=0,\quad L_{-}yQ=-2\nabla Q, \quad L_{-}|y|^{2}Q=-4\Lambda Q.	
\end{equation}
Here $\rho$ is Schwarz and decays exponentially fast.\\
Note that the last formula of \eqref{eq: spectralba1} and \eqref{eq: spectralba2} are directly related to the pseudo-conformal symmetry of \eqref{eq: nls}.\\
\subsubsection{Coersiveness}
We recall here the important coersiveness property for $L_{+}, L_{-}$, we refer to \cite{merle2005blow},\cite{merle2006sharp},\cite{weinstein1986lyapunov}, and see Lemma 3.2 in \cite{raphael2011existence}.
\begin{lemma}[Lemma 3.2 in \cite{raphael2011existence}]\label{lem: coer}
There exists $\delta>0,C>0$, such that for all $\epsilon:=\epsilon_{1}+i\epsilon_{2}\in H^{1}$, one has 
\begin{equation}	
\begin{aligned}
&<L_{+}\epsilon_{1}, \epsilon_{1}>+<L_{-}\epsilon_{2},\epsilon_{2}>\\
\geq &\delta \|\epsilon\|_{H^{1}}^{2}-C\left(<\epsilon_{1},Q>^{2}+<\epsilon_{1}, |y|^{2}Q>^{2}+<\epsilon_{1}, yQ>^{2}+<\epsilon_{2},Q>^{2}\right)
\end{aligned}
\end{equation}
\end{lemma}

\subsection{The elliptic object}
Let $\PPP=(\lambda,b,\alpha,\beta)$ be parameters, where $\lambda\in \mathbb{R+}, b\in \mathbb{R}, \alpha,\beta\in \mathbb{R}^{2}$. There exists a family of elliptic objects $Q_{\PPP}$, which are close to $Q$, and play a crucial role in the construction of minimal blow up of \cite{raphael2011existence}. We recall its property in this subsection. \\

It should be noted, in \cite{raphael2011existence}, because $k\in C^{5}$, this allows the authors to construct a high order (up to fifth order) approximation of $Q_{\PPP}$ to $Q$.
 Here, since our $k$ is only $C^{2}$, the elliptic object $Q_{\PPP}$ we use this article is indeed only includes the up to second order expansion in the construction of \cite{raphael2011existence}, i.e. the $Q_{\PPP}$
 in this article is not exactly the same $Q_{\PPP}$ as in \cite{raphael2011existence}. In some sense, our  $Q_{\PPP}$ is simpler than the one in \cite{raphael2011existence}, thus simplified part of the computations in \cite{raphael2011existence}, but this has to been compensated by further modifying the Lyapounov control in \cite{raphael2011existence}.
 
We use $|P|:=|\lambda|+|b|+|\alpha|+|\beta|$ to measure the size of those modulation parameters. We summarize the basic properties of $Q_{\PPP}$ as,

\begin{lemma}\label{lem: qp}[Prop 2.1 and section 2.2 in  \cite{raphael2011existence}]

For $|\PPP|$ small, there exists $Q_{\PPP}$ with, 
\begin{equation}\label{eq: qp}
\begin{aligned}
Q_{\PPP}=P_{\PPP}e^{-ib\frac{|y|^{2}}{4}+i\beta y},\\
P_{\PPP}=Q+T_{2}\\
\end{aligned}
\end{equation}
where $T_{2}$ satisfying 
\begin{equation}\label{eq: t2}
\begin{aligned}
L_{+}T_{2}=\nabla^{2}k(0)(\alpha,y)\lambda Q^{3}+\frac{\lambda^{2}}{2}\nabla^{2}k(0)(y,y)Q^{3}-\lambda c_{0}(\alpha)yQ
\end{aligned}
\end{equation}
where 
\begin{equation}\label{eq: coa}
c_{0}(\alpha)_{j}=\frac{\int Q^{4}}{2\int Q^{2}}\nabla^{2}k(0)(e_{j},\alpha), \quad j=1,2.
\end{equation}
Furthermore, $Q_{\PPP}$ solves
\begin{equation}\label{eq: qb}
\begin{aligned}
-&b^{2}\partial_{b}Q_{\PPP}+i(-b\beta+\lambda c_{0}(\alpha))\partial_{\beta}Q_{\PPP}\\
+&\Delta Q_{\PPP}-Q_{\PPP}+\frac{k(\lambda y+\alpha)}{k(\alpha)}|Q_{\PPP}|^{2}Q_{\PPP}+ib\Lambda Q_{\PPP}-2i\beta \nabla Q_{\PPP}-|\beta|^{2}Q_{P}\\
=-&\Psi_{\PPP}
\end{aligned}
\end{equation}

with estimates (for some constant $c>0$)
\begin{equation}\label{eq: estimateforqb}
|Q_{\PPP}|+|\partial Q_{\PPP}|+|\partial^{2}Q_{\PPP}|\lesssim e^{-c|y|} 
\end{equation}
and\footnote{It should be noted that $\Psi$ and $\tilde{\Psi}$  satisfy the same type of estimates. }
\begin{equation}\label{eq: estimateforphib}
|\tilde{\Psi}_{\PPP}|+|\partial\tilde{\Psi}_{\PPP}|+|\partial^{2}\tilde{\Psi}_{\PPP}|\lesssim o_{\PPP}(\PPP^{2})e^{-c|y|} 
\end{equation}
where $\tilde{\Psi}_{\PPP}:=e^{i\frac{b|y|^{2}}{4}-i\beta y}\Psi_{\PPP}$.

And as observed in \cite{raphael2011existence}, one has 
\begin{equation}\label{eq: algeor}
	<T_{2},Q>=0.
\end{equation}
\end{lemma}

\begin{remark}
Strictly speaking, it may be better to denote $T_{2}$ as $T_{2,\PPP}$ since it does depend on the parameter $\PPP$, but we will follow the notation in \cite{raphael2011existence}.

 The equation for $T_{2}$ is (2.12) in \cite{raphael2011existence} and condition for $c_{0}(\alpha)$ is to ensure the existence of $T_{2}$ by Fredholm alternative $<L_{+}T_{2}, \nabla Q>=0.$\\
 
 The equation \eqref{eq: estimateforphib} is  (2.22) in 
 \cite{raphael2011existence}, but we absorb $O(\PPP^{3})$ terms into $\Psi_{b}$. It can be (algebraically) checked that, (up to $O(\PPP^{3})$), \eqref{eq: estimateforphib} is equivalent to
 \begin{equation}\label{eq: pb}
 \begin{aligned}
 -&\lambda c_{0}(\alpha)yP_{\PPP}\\
 +&\Delta P_{\PPP}-P_{\PPP}+\frac{k(\lambda y+\alpha)}{k(\alpha)}|P_{\PPP}|^{2}P_{\PPP}\\
 =&-\tilde{\Psi}_{\PPP}(y)+O(P^{3})e^{-c|y|}
 \end{aligned}
 \end{equation}
Note that \eqref{eq: pb} is (up to $O(P^{3})$ error), (2.8) in \cite{raphael2011existence}.

\end{remark}

\begin{remark}
Note the parameters $\lambda, b, \alpha,\beta$ are corresponding to scaling, pesudo-conformal symmetry, space translation, and galiean symmetry. In some sense, the construction of $Q_{\PPP} $ is remedifying the loss of those exact symmetries for \eqref{eq: knls}.
\end{remark}

We finally recall the mass and energy satisfy
\begin{lemma}\label{lem: me}[Analogue of Lemma 2.4 in \cite{raphael2011existence}]
\begin{enumerate}
	\item Mass
	\begin{equation}\label{eq: qbmass}
	\int |Q_{\PPP}|^{2}=\int Q^{2}+O(P^{4})
	\end{equation}
	\item Energy
	\begin{equation}\label{eq: qbenergy}
	\tilde{E}(Q_{\PPP}):=\frac{1}{2}\int |\nabla Q_{\PPP}|^{2}-\frac{1}{4}\int \frac{k(\lambda y+\alpha)}{k(\alpha)}|Q_{\PPP}|^{4}=
	\frac{b^{2}}{8}\int |y|^{2}Q^{2}+\frac{\beta^{2}}{2}\int Q^{2}-\frac{\lambda^{2}}{8}\int \nabla^{2}k(0)(y,y)Q^{4}+o_{\PPP}(P^{2}).
	\end{equation}
\end{enumerate}	
\end{lemma}
We note since $Q_{\PPP}$ in the current article only matches the $Q_{\PPP}$ up to $O(\PPP^{3})$, \eqref{eq: qbenergy} is slightly different to the one in Lemma 2.4 in \cite{raphael2011existence}.

Lemma \ref{lem: me} follows mainly from $<T_{2},Q>=0$ and $-\Delta Q+Q=Q^{3}$, and the little o notation $o_{\PPP}(P^{2})$ follows from $k(z)=k(0)+\nabla^{2}k(0)(z,z)+o_{z}(|z|^{2})$. We refer to the proof of Lemma 2.4 in \cite{raphael2011existence} for more details.

\subsection{The modulation equations}
In this subsection, we recall the modulation equation for ansatz
\begin{equation}\label{eq: ansatzmain}
\begin{aligned}
u(t,x)&:=\frac{1}{k^{\frac{1}{2}}(\alpha)}\frac{1}{\lambda(t)}(Q_{\PPP}+\epsilon)(y)e^{i\gamma(t)},\\
\frac{ds}{dt}&=\frac{1}{\lambda^{2}}, \quad y=\frac{x-\alpha(t)}{\lambda(t)}
\end{aligned}
\end{equation}

It should be noted, since our $Q_{\PPP}$ is not exactly \footnote{Recall due to the fact our $k(x)$ is only $C^{2}$, one cannot approximate $Q$ as in \cite{raphael2011existence} up to $O(\PPP^{5})$ error.}the same $Q_{\PPP}$ as in \cite{raphael2011existence}, strictly speaking, both ansatz \eqref{eq: ansatzmain} and the associated modulational equations are not exactly the same one as \cite{raphael2011existence}, (and they admits different estimates), nevertheless, a lot of algebraic computation are exactly same, and formally, the associated modulation equations are same. 

Viewing \eqref{eq: knls} as the evolution of $\epsilon$ coupled with the evolution of parameters\footnote{We will, as most articles in this field, abuse notation a bit to write, for example, $b(t(s))$ simply as $b(s)$.} $\lambda(t), b(t),\alpha(t),
\beta(t)$, one has 

\begin{equation}\label{eq: ep1}
\begin{aligned}	
&i\partial_{s}\epsilon+\Delta \epsilon-\epsilon+\frac{k(\lambda y+\alpha)}{k(\alpha)}(|Q_{\PPP}+\epsilon|^{2}(Q_{\PPP}+\epsilon)-|Q_{\PPP}|^{2}Q_{\PPP})\\
+&ib\Lambda \epsilon-i2\beta(\nabla \epsilon+\frac{\lambda}{2}\frac{\nabla k(\alpha)}{k(\alpha)}\epsilon)-|\beta|^{2}\epsilon\\
=&-i(\frac{\lambda_{s}}{\lambda}+b)\lambda\partial_{\lambda}Q_{\PPP}-i(b_{s}+b^{2})\partial_{b}Q_{\PPP}-i(\frac{\alpha_{s}}{\lambda}-2\beta)\lambda\partial_{\alpha}Q_{\PPP}-i(\beta_{s}+b\beta-c_{0}(\alpha)\lambda)\partial_{\beta}Q_{\PPP}\\
+ & (\tilde{\gamma}_{s}-|\beta|^{2})(Q_{\PPP}+\epsilon)+i(\frac{\lambda_{s}}{\lambda}+b)(\Lambda Q_{\PPP}+\Lambda\epsilon)+i(\frac{\alpha_{s}}{\lambda}-2\beta)[\nabla(Q_{\PPP}+\epsilon)+\frac{\lambda}{2}\frac{\nabla k(\alpha)}{k(\alpha)}(Q_{\PPP}+\epsilon)] \\
+ & \Psi_{\PPP}+O(\PPP^{3})e^{-c|y|},
\end{aligned}
\end{equation}
Note that if one takes the Real part and Imaginary part of \eqref{eq: ep1}, it is (up to $O(P^{3})$ error) (3.11) and (3.12) in \cite{raphael2011existence}.

We finally recall the orthogonality condition in \cite{raphael2011existence}. It follows from implicit function theorem, (up to checking some standard non-degenerate condition), 
\begin{lemma}\label{lem: impli}
There exists $\delta_{0}>0$, so that for any $f\in H^{1}$ with
\begin{equation}
f=\frac{1}{k^{1/2}(\tilde{\alpha})}\frac{1}{\tilde{\lambda}}(Q_{\tilde{\PPP}}+\tilde{\epsilon})(\frac{x-\tilde{\alpha}}{\tilde{\lambda}})e^{i\tilde{\gamma}}
\end{equation}

with
\begin{equation}\label{eq: initialsmall}
|\tilde{\PPP}|+\|\tilde{\epsilon}\|_{H^{1}}\leq \delta\leq \delta_{0},	
\end{equation}
then there exist unique $\tilde{\PPP}:=(\tilde{\lambda}, \tilde{b},\tilde{\alpha}, \tilde{\beta})$, and $\tilde{epsilon}\in H^{1}$, $\tilde{\gamma}$, with
\begin{enumerate}
\item 	Smallness condition
\begin{equation}\label{eq: smallness}
|\PPP|+\|\epsilon\|_{H^{1}}
+|\gamma-\tilde{\gamma}|\lesssim \delta.
\end{equation}
\item  Decomposition
\begin{equation}
	f=\frac{1}{k^{1/2}(\alpha)}\frac{1}{\lambda}(Q_{\PPP}+\epsilon)(\frac{x-\alpha}{\lambda})e^{i\gamma}
\end{equation}
\item Orthogonality condition
\begin{equation}
\Im <\nabla Q_{\PPP},\epsilon>=
(\epsilon_{2},\nabla\Sigma)-(\epsilon_{1},\nabla\Theta)=0,\label{eq: orth1}
\end{equation}
\begin{equation}
\Re <yQ_{\PPP},\epsilon>=
(\epsilon_{1},y\Sigma)+(\epsilon_{2},y\Theta)=0,\label{eq: orth2}
\end{equation}
\begin{equation}
\Im<\Lambda Q_{\PPP},\epsilon>=
-(\epsilon_{1},\Lambda\Theta)+(\epsilon_{2},\Lambda\Sigma)=0,\label{eq: orth3}
\end{equation}
\begin{equation}
\Re<y^{2}Q_{\PPP},\epsilon>=
(\epsilon_{1},|y|^{2}\Sigma)+(\epsilon_{2},|y|^{2}\Theta)=0,\label{eq: orth4}
\end{equation}
\begin{equation}
\Im<\epsilon, \rho e^{-i\frac{b|y|^{2}}{4}+i\beta y}>=
-(\epsilon_{1},\rho_{2})+(\epsilon_{2},\rho_{1})=0,\label{eq: orth5}
\end{equation}
\end{enumerate}

\end{lemma}
Here we use notation $Q_{\PPP}:=\Sigma+i\Theta$, and $\rho e^{-i\frac{b|y|^{2}}{4}+i\beta y}=\rho_{1}+i\rho_{2}$.\\
One may refer to \cite{merle2005blow},\cite{raphael2011existence} for a proof. Note that Lemma \ref{lem: impli} is purely stationary, with no evolution in time involved.

On the other hand, via local wellposedness of \eqref{eq: knls} in $H^{1}$, one has 
\begin{lemma}\label{lem: lwpimpli}
Let $\delta_{0}$ be as in Lemma \ref{lem: impli}. Let $u$ solves \eqref{eq: knls} with initial data $u_{0}$ such that 
\begin{equation}\label{eq: geod}
u_{0}=\frac{1}{k^{1/2}(\alpha_{0})}\frac{1}{\lambda_{0}}(Q_{\PPP_{0}}+\epsilon_{0})(\frac{x-\alpha_{0}}{\lambda_{0}})e^{i\gamma_{0}}	
\end{equation}
with smallness condition \eqref{eq: initialsmall} and orthogonality condition \eqref{eq: orth1}-\eqref{eq: orth5}. Then there exists $T=T(\|u_{0}\|_{H^{1}})>0$, so that for all $t\in [0,T]$, one has a unique decomposition as 
\begin{equation}\label{eq: gforawhile}
u(t,x)=\frac{1}{k^{1/2}(\alpha(t))}\frac{1}{\lambda(t)}(Q_{\PPP(t)}+\epsilon)(\frac{x-\alpha(t)}{\lambda(t)})e^{i\gamma(t)}
\end{equation}
with \eqref{eq: smallness}, and orthogonality condition \eqref{eq: orth1}-\eqref{eq: orth5}.
\end{lemma}

\subsection{Motivation for the construction of $Q_{\PPP}$}
We briefly explain the motivation of the construction of $Q_{\PPP}$, we refer to \cite{raphael2011existence} for more details.
Consider ansatz \eqref{eq: ansatzmain} with $\epsilon\equiv 0$,
and following the modulation equation for the $S(t)$ blow up for \eqref{eq: nls},
\begin{equation}\label{eq: mo}
\begin{cases}
\frac{\lambda_{s}}{\lambda}+b=0,\\
\frac{\alpha_{s}}{\lambda}-2\beta=0,\\
\tilde{\gamma}_{s}-|\beta|^{2}=0, & \tilde{\gamma}_{s}:=\gamma_{s}-1,
\end{cases}
\end{equation}

One obtains \eqref{eq: qb} up to correction terms. The $c_{0}(\alpha)$ term in \eqref{eq: qb}, as aforementioned, is a technical correct to admit certain Fredhlom alternative.

In some sense, the construction of $Q_{\PPP}$ is to ensure \eqref{eq: mo} holds for the whole ansatz \eqref{eq: ansatzmain}, at least approximately up to correction terms.

\section{Setting up and main ansatz}\label{sec: mainans}
\subsection{Overview}
As aforementioned, the overall strategy, following \cite{merle1990construction}, \cite{raphael2011existence},  to construct the blow up solution in Theorem \ref{thm: main} is to solve a sequence of solution $u_{n}$ to \eqref{eq: knls} backwards with $u_{n}(t_{n})=u_{0,n}$, with $t_{n}\rightarrow 0$, and $u_{0,n}$ well prescribed, so that one can obtain uniform estimates for $u_{n}(T)$ for some $T<0$, and establish some compactness so that (up to extracting subsequence), $u_{n}(T)$ converges to $u_{0,c}$, and solution $u_{c}$ to \eqref{eq: knls} with $u_{c}(T):=u_{0,c}$ will give a desired blow up solutions to \eqref{eq: knls}. This procedure is more or less standard, giving one can perform energy type arguments detailed dynamic of $u_{n}$. We will sketch the compactness argument in Subsection \ref{sub: compactness}.  We focus on description and $u_{0,n}$ and the associated bootstrap Lemma describing its dynamic.

\subsection{Bootstrap Setting up}\label{subsec: boot}
Let $t_{n}\rightarrow 0$, $n\geq 1$, and we aim to study solution $u_{n}$ to \eqref{eq: knls} with prepared initial data $u_{n}(t_{n})=u_{0,n}$ within $[T_{0}, t_{n}]$ for some $T_{0}<0$, and get some uniform estimates for $u_{n}(T_{0})$, (in dependent of $n$). Before we go to more details ( and more specific description), for notation convenience, we will fix $n$, and set $t_{n}=t_{0}$, and we will drop the $n$ notation in $u_{n}, u_{0,n}$, and simply denote it as $u, u_{0}$.

We fix a small parameter $\delta$ thorough the whole section. And for bootstrap purpose, we also fix a large parameter $K$ throughout the whole section.
\subsubsection{Description of initial data}\label{subsub: initial}
Let $u_{0}$ be of form \eqref{eq: geod}. And we will pose $\epsilon_{0}\equiv 0$, $\alpha_{0},\beta_{0}=0$,\footnote{There is some freedom to choose the initial data. We choose zero initial data just for convenience. }  and 
\begin{equation}\label{eq: initialdata}
u_{0}=\frac{1}{\lambda_{0}}(Q_{\PPP_{0}})(\frac{x-\alpha_{0}}{\lambda_{0}})e^{i\gamma_{0}}	
\end{equation}
Thus orthogonality conditions \eqref{eq: orth1}-\eqref{eq: orth5} automatically holds.

We always assume
\begin{equation}\label{eq: smallnessini}
|\PPP|\ll 1, |t_{0}|\ll 1.	
\end{equation}

It should be noted one has , by Lemma \ref{lem: me}, one has\footnote{Strictly speaking, we should bound LHS of \eqref{eq: massini} by $\PPP_{0}^{4}$, but since we will also assume later \eqref{eq: inimore}, we can bound $\PPP_{0}$ via $\lambda_{0}$. } 
\begin{equation}\label{eq: massini}
|\|u_{0}\|_{L_{x}^{2}}^{2}-\|Q\|_{L_{x}^{2}}|\lesssim \lambda_{0}^{4}	
\end{equation}

We further assume
\begin{equation}\label{eq: inimore}
b_{0}=-\frac{t_{0}}{C_{0}^{2}}, \lambda_{0}=-\frac{t_{0}}{C_{0}}, \alpha_{0}=\beta_{0}=0, \gamma_{0}=-\frac{C_{0}^{2}}{t_{0}}.
\end{equation}
Note that, we may also denote $b_{0}=b(t_{0}), \lambda_{0}=\lambda(t_{0})$,etc.

Also note that, \eqref{eq: inimore} in particular implies
\begin{equation}\label{eq: estimateini}
|\frac{\beta_{0}}{\lambda_{0}}|+|\frac{\alpha_{0}}{\lambda_{0}}|+|\lambda_{0}+\frac{t_{0}}{C_{0}}|+|\frac{b_{0}}{\lambda_{0}}-\frac{1}{C_{0}}|=0.
\end{equation}

Here we have the freedom to prescribe $C_{0}>0$.

We finally note that, \eqref{eq: inimore} gives $\lambda_{0}=C_{0}b_{0}$
Furthermore, by \eqref{eq: qbenergy}, and \eqref{eq: smallnessini}, we have
\begin{equation}\label{eq: energycondition}
\tilde{E}_{0}:=E_{in}(u_{0})+\frac{1}{8}\int \nabla^{2}k(0)(y,y)Q^{4}=\frac{b_{0}^{2}}{\lambda_{0}^{2}}\
||yQ\|_{2}^{2}+o_{\PPP}(\mathcal{P}^{2})>0
\end{equation}
Note that \eqref{eq: energycondition} is crucial in \cite{raphael2011existence}.
\begin{remark}
The energy shift phenomenon is caused by the non-degenerate condition in \eqref{as: main}. It is also a necessary condition for the existence of a critical blow-up solution. One can refer to Prop 1.2 and Comment 1 for Theorem 1.3 in \cite{raphael2011existence}.
\end{remark}

And we note
\begin{equation}\label{eq: atecnicalconstant}
C_{0}=\frac{\|yQ\|_{L^{2}}}{\sqrt{8\tilde{E}_{0}}}+o_{\PPP}(1)
\end{equation}
\begin{remark}\label{rem: extendconst}
We note that the introduction of this constant $\frac{\|yQ\|_{L^{2}}}{\sqrt{8\tilde{E}_{0}}}$, (as well as $\tilde{E}_{0}$) is very natural in the sense since by \eqref{eq: conesti}, one formally has, (up to lower order terms),
\begin{equation}
b^{2}\|yQ\|_{2}^{2}\approx \lambda^{2}(E(u_{0})+\frac{1}{8}\int\nabla^{2}k(0)(y,y)Q^{4}).
\end{equation} 
\end{remark}

\subsubsection{Main bootstrap Lemma}
We are now ready to state our main bootstrap Lemma. Let $u$ be associated solution to \eqref{eq: knls} with $u(t_{0})=u_{0}$. By Lemma \ref{lem: lwpimpli}, we know at least for $T$ close to $t_{0}$, one has unique geometric decomposition for $u$ in $[T,t_{0}]$ as \eqref{eq: gforawhile}, with orthogonality condition \eqref{eq: orth1}-\eqref{eq: orth5}.
\begin{equation}
u(t,x)=\frac{1}{k^{1/2}(\alpha(t))}\frac{1}{\lambda(t)}(Q_{\PPP(t)}+\epsilon)(\frac{x-\alpha(t)}{\lambda(t)})e^{i\gamma(t)}	
\end{equation}

We have
\begin{lemma}\label{lem: bootmain}
Let $u$ be associated solution to \eqref{eq: knls} with $u(t_{0})=u_{0}$, and $\delta$ be a small parameter. We assume $|T|\ll 1$.

Assume bootstrap assumption for $u(t)$, $t$ in $[T, t_{0}]$, 
\begin{equation}\label{eq: baep}
\|\epsilon(t)\|_{H^{1}}\leq \lambda(t),
\end{equation}
\begin{equation}\label{eq: bamome}
|\frac{\beta(t)}{\lambda(t)}|+|\frac{\alpha(t)}{\lambda(t)}|\leq \delta,	
\end{equation}
\begin{equation}\label{eq: balam}
|\lambda(t)+\frac{t}{C_{0}}|\leq \delta \lambda(t)
\end{equation}
\begin{equation}\label{eq: bab}	
|\frac{b(t)}{\lambda(t)}-\frac{1}{C_{0}}|\leq \delta.
\end{equation}

Then, one has bootstrap estimates
\begin{equation}\label{eq: beep}
\|\epsilon(t)\|_{H^{1}}\leq \lambda^{\frac{3}{2}}(t),
\end{equation}
\begin{equation}\label{eq: bemome}
|\frac{\beta(t)}{\lambda(t)}|+|\frac{\alpha(t)}{\lambda(t)}|\leq \frac{1}{2}\delta,	
\end{equation}
\begin{equation}\label{eq: belam}
|\lambda(t)+\frac{t}{C_{0}}|\leq \frac{1}{2}\delta \lambda(t)
\end{equation}
\begin{equation}\label{eq: beb}	
|\frac{b(t)}{\lambda(t)}-\frac{1}{C_{0}}|\leq \frac{1}{2}\delta.
\end{equation}
\end{lemma}

We note that the bootstrap assumptions implies, in particular, for all $t\in [T,t_{0}]$, one has 
\begin{equation}\label{eq: preboot1}
|\PPP(t)|\lesssim \lambda(t)\sim b(t)	
\end{equation}
and 
\begin{equation}\label{eq: preboot2}
\lambda(0)\lesssim \lambda(t), \quad b(0)\lesssim b(t).	
\end{equation}

\begin{remark}
By the standard continuity argument, one can keep evolving $u$ backwards as far as bootstrap estimates ensures $\|\PPP\|+\|\epsilon\|_{H^{1}}\ll 1$ so that Lemma \ref{lem: impli} works, this will give the  $T_{0}$ (independent of $t_{0}$ as far as it is close to zero), such that $u(T_{0})$ admits uniform estimates, as mentioned in the beginning of Subsection \ref{subsec: boot}.

We note that estimate \eqref{eq: estimateini} and the fact $\epsilon_{0}=0$, ensures bootstrap estimates always hold for $T$ close to $t_{0}$ enough (by local wellposedness)
\end{remark}
\subsection{A sketch of compactness argument}\label{sub: compactness}
Let $t_{n}$ $\rightarrow$ $0$ and $u_{0,n}:=\frac{1}{\lambda_{n}}(Q_{\PPP_{n}}(\frac{x}{\lambda_{n}}))e^{i\gamma_{n}}$ satisfying the description of  initial data in Subsection \ref{sec: bootstrap}, so that Lemma \ref{lem: bootmain} holds. Then let $T_{0}<0$, so that Lemma \ref{lem: bootmain} works within $[T_{0}, t_{n}]$. It follows from Lemma \ref{lem: bootmain} that $u_{n}(T_{0})$ are uniform bounded in $H^{1}$. Note that $H^{1}\rightarrow L^{2}$ is not compact due to space translation. However, as observed in \cite{merle1990construction}, \cite{raphael2011existence}, one may study the local mass $\int \chi(R)|u_{n}(T_{0})|^{2}$, $R\geq 1$.

Indeed, one has 
\begin{equation}\label{eq: iniextmass}
\int \chi(R)|u_{n}(t_{n})|^{2}\rightarrow 0,	
\end{equation}

And 
\begin{equation}\label{eq: slowchange}
\partial_{t}\int \chi(R)|u_{n}(t_{n})|^{2}=\Im\int \frac{1}{R}\nabla \chi(\frac{x}{R})\nabla u_{n}\bar{u}_{n}=O(\frac{1}{R})
\end{equation}

Note that \eqref{eq: slowchange} does not follow from $u_{n}(t)$, $t\in [T_{0},t_{n}]$, are uniform bounded in $H^{1}$, and they are \textbf{not}. But from Lemma \ref{lem: bootmain}, we know $u_{n}(t)-\frac{1}{\lambda_{n}(t)}Q(x/\lambda(t))e^{i\gamma_{n}(t)}$ is bounded in $H^{1}$, since $Q$ is fast decaying, thus \eqref{eq: slowchange}.

Since $u_{n}(T_{0})$ are uniform bounded in $H^{1}$, and \eqref{eq: iniextmass} and \eqref{eq: slowchange} implies 
the mass cannot escape to infinite. Up to picking subsequence, we know $u_{n}(T_{0})$ converges to some $u_{0,c}$ in $L^{2}$ with $\|u_{0,c}\|_{2}=\|Q\|_{2}$. And by the $L^{2}$ LWP, one obtains that $u_{c}$ blows up at $T=0$.

With the above choice of approximate profile $Q_{\PPP}$ and orthogonality conditions, we are ready to give a bootstrap argument in this section. The aim is to restrict the solution to evolve in the direction as desired. Once we finish the bootstrap argument, the standard compactness argument helps us construct the minimal mass blow up solution in Theorem \ref{thm: main}.

\begin{remark}
Further concentration compactness can yield that $u_{0,c}$ is of blow up rate $\sim \frac{1}{t}$.
\end{remark}

\section{Proof of Bootstrap Lemma}\label{sec: bootstrap}
Let $u$ be as in Lemma \ref{lem: bootmain}, and we assume $t\in [T,t_{0}]$, thus all bootstrap assumptions \eqref{eq: baep}-\eqref{eq: bab} holds.

In this section, we first collect estimates induced by conservation laws and modulational estimates, which comes from the choice of (non-degenerate) orthogonality conditions. Then we prove Lyapounov type energy estimates, which are the key to control $\|\epsilon\|_{H^{1}}$. Finally, we prove bootstrap estimates in Lemma \ref{lem: bootmain} and fill in other technical details.

\subsection{Estimates by conservation laws and modulational estimates}

Firstly, by energy conservation and mass conservation law, one has 
\begin{lemma}\label{lem: conesti}
Let $u$ be as in Lemma \ref{lem: bootmain}, one has 
\begin{equation}\label{eq: conesti}
  \frac{b^{2}}{8}\int|y|^{2}Q^{2}+\frac{|\beta|^{2}}{2}\int Q^{2}-\frac{\nabla^{2}k(0)(\alpha,\alpha)}{4}\int Q^{2}-\lambda^{2}(E_{0}+\frac{1}{8}\int\nabla^{2}k(0)(y,y)Q^{4})= \|\epsilon\|_{H^{1}}^{2}+O(\PPP^{2}\|\epsilon\|_{H_{1}})+o_{\PPP}(\PPP^{2})	
\end{equation}
\end{lemma}

Secondly, by using modulational equation \eqref{eq: ep1} and plug in orthogonality conditions \eqref{eq: orth1}-\eqref{eq: orth5}, one has 
\begin{lemma}\label{lem: moduesti}
Let $u$ be as in Lemma \ref{lem: bootmain}, one has 

\begin{equation}\label{eq: moduesti}
\begin{aligned}
&|b_{s}+b^{2}-d_{0}(\alpha,\alpha)|+|\tilde{\gamma}_{s}-|\beta|^{2}|+|\frac{\alpha_{s}}{\lambda}-2\beta|+|\frac{\lambda_{s}}{\lambda}+b|+|\beta_{s}+b\beta-c_{0}(\alpha)\lambda|\\
\leq &\|\epsilon\|_{H^{1}}^{2}+\PPP^{2}\|\epsilon\|_{H^{1}}+o_{\PPP}(1)\PPP^{2}.
\end{aligned}	
\end{equation}
In particular, one has the following crude bound
\begin{equation}\label{eq: crude}
\partial_{s}\PPP=O(\PPP^{2}).
\end{equation}
\end{lemma}

Here, as in \cite{raphael2011existence}, we define
\begin{equation}
	d_{0}(\alpha,\alpha)=\frac{2||Q||_{L^{2}}^{2}}{||yQ||_{L^{2}}^{2}}\nabla^{2}k(0)(\alpha,\alpha),\quad d_{1}(\alpha,\alpha)=\frac{(|y|^{2}Q,\rho)}{4(\rho,Q)}d_{0}(\alpha,\alpha).
\end{equation}

Lemma \ref{lem: conesti}, \ref{lem: moduesti} are the analogue of Lemma 3.1 in \cite{raphael2011existence}.

We will prove Lemma \ref{lem: conesti}, \ref{lem: moduesti}
in Subsection \ref{subsec: tech}.

\subsection{Lyapounov type energy estimates}\label{subsec: lya}
\subsubsection{Overview}
The goal of the current subsection is to construct a functional $I_{1}(\epsilon)$, such that 

\begin{lemma}\label{lem: keyly}
For all $t\in [T,t_{0}]$, under the assumptions of Lemma \ref{lem: bootmain}, one has 
\begin{enumerate}
	\item $I_{1}$ is almost non-decreasing in the sense\footnote{Here, we say $A\geq O(1)$ if there exists $C>0$ so that $A\geq -C.$}
	\begin{equation}
	\frac{d}{dt}I_{1}(\epsilon)\equiv \frac{1}{\lambda^{2}}\frac{d}{ds}I_{1}(\epsilon)\geq O(1).
	\end{equation}
	\item $I_{1}$ is almost coercive in the sense, there exists $\delta_{0}>0$,
	\begin{equation}\label{eq: lc3}
	\lambda^{2}I_{1}(\epsilon)\geq \delta_{0}\|\epsilon\|_{H^{1}}^{2}+O(\mathcal \PPP^{3}).	
	\end{equation}
\end{enumerate}

\end{lemma}

$I_{1}$ will be constructed based on a modification of Lyapounov functional in \cite{raphael2011existence}.\\
 
Such modifications are necessary to handle $k\in C^{2}$ mainly because in the current (critically) low regularity $k$ situation, it is impossible to construct $Q_{\PPP}$ approximate to $Q$ to high order as in \cite{raphael2011existence}, and it in turn makes $\Psi_{\PPP}$ large.

Recall anstatz \eqref{eq: ansatzmain}, in \cite{raphael2011existence}, the Lyapounov functional was based on $\tilde{u}:=\frac{1}{\lambda}\epsilon(\frac{x-\alpha}{\lambda})e^{i\gamma}$. We will formulate the Lyapounov functional based on $\epsilon$. On one hand, this  gives another view point for the Virial algebra used in \cite{raphael2011existence}, which may be of interest. On the other hand, it is more natural for us to present the further modification based on this $\epsilon$ formulation. Before we go to more details, we mention that when we say Virial algebra, we mean $\int \Lambda f 
\Delta \bar{f} $ has a sign for all $f$.

\subsubsection{A further simplification of \eqref{eq: ep1} with a revisit of the Lyapounov functional in \cite{raphael2005stability}}
We first simplify \eqref{eq: ep1} a bit. Morally, as one will see from the computations below, any term of $O(\PPP^{3})$ should be understood as purely perturbative in \eqref{eq: ep1}, and by bootstrap assumption \eqref{eq: baep}, terms of $O(P^{2}\epsilon)$, etc should also be viewed of purely petrubative. 

Recall Lemma \ref{lem: qp}, we rewrite \eqref{eq: ep1} as 
\begin{equation}\label{eq: ep1s1}
\begin{aligned}
&i\partial_{s}\epsilon+\Delta \epsilon-\epsilon+(|Q_{\PPP}+\epsilon|^{2}(Q_{\PPP}+\epsilon))-|Q_{\PPP}|^{2}Q_{\PPP})+ib\Lambda \epsilon-i2\beta\nabla \epsilon\\
=&i(b_{s}+b^{2})\frac{y^{2}Q_{\PPP}}{4}-i(\beta_{s}+b\beta-c_{0}(\alpha)\lambda)yQ_{\PPP}+(\tilde{\gamma}_{s}-\beta^{2})Q_{\PPP}+i(\frac{\lambda_{s}}{\lambda}+b)(\Lambda Q_{\PPP}+\epsilon)\\+&i(\frac{\alpha_{s}}{\lambda}-2\beta)\nabla( Q_{\PPP}+\epsilon)\\
+&\Psi_{\PPP}+F+(1-\frac{k((\lambda y+\alpha)}{k(\alpha)})|\epsilon|^{2}\epsilon
\end{aligned}
\end{equation}
where we denote error terms as $F$,
\begin{equation}\label{eq: F}
	F=O(\PPP^{3})e^{-c|y|}+O(\PPP^{2}
	\epsilon).
\end{equation}
\\
We further simplify \eqref{eq: ep1s1} as
\begin{equation}\label{eq: ep1s2}
\begin{aligned}
&i\partial_{s}\epsilon+\Delta \epsilon-\epsilon+(|Q_{\PPP}+\epsilon|^{2}(Q_{\PPP}+\epsilon))-|Q_{\PPP}|^{2}Q_{\PPP})+ib\Lambda \epsilon-i2\beta\nabla \epsilon\\
=&O(\PPP^{2})(y^{2}Q_{\PPP}+yQ_{\PPP}+Q_{\PPP}+i\Lambda Q_{\PPP}+i\nabla Q_{\PPP})\\
+&O(\PPP^{2}i\Lambda \epsilon+O(\PPP^{2})i\nabla \epsilon\\
+&\Psi_{\PPP}+F+(1-\frac{k(\lambda y+\alpha)}{k(\alpha)})|\epsilon|^{2}\epsilon.
\end{aligned}
\end{equation}

Let us denote the term first order in $\epsilon$ of $-\Delta \epsilon+\epsilon-(|Q_{\PPP}+\epsilon|^{2}(Q_{\PPP}+\epsilon))-|Q_{\PPP}|^{2}Q_{\PPP}$ as $M(\epsilon)$, then one has 
\begin{equation}
\begin{aligned}
M(\epsilon)=M_{1}(\epsilon)+iM_{2}(\epsilon),
\end{aligned}	
\end{equation}
This operator $M$ is self adjoint (with respect to inner product $\Re <f,g>:\equiv \Re \int f\bar{g}$) is the perturbed version of the Sch\"odinger operator $L(\epsilon):=L_{+}\epsilon_{1}+iL_{-}\epsilon_{2}$, and has been well studied in \cite{raphael2011existence}. We note that it may be better to denote $M$ as $M_{\PPP}$, but we follow the notation conventions in the literature omitting this $\PPP$.

We also denote the higher order in $\epsilon$ as $R(\epsilon)$, i.e
\begin{equation}
-\Delta \epsilon+\epsilon-(|Q_{\PPP}+\epsilon|^{2}(Q_{\PPP}+\epsilon))-|Q_{\PPP}|^{2}Q_{\PPP}=M(\epsilon)+R(\epsilon).	
\end{equation}

\begin{remark}
 Since we mainly want to focus a revisit of Lyapounov functional in \cite{raphael2011existence}, and highlight the difference of our version and further modifications, we are oversimplifying a bit in \eqref{eq: ep1s2}. Strictly speaking, one needs the $O(\PPP^{2})i\Lambda \epsilon$ to be exactly $(b+\lambda_{s}/\lambda)$ in \eqref{eq: ep1s2}, (and $O(\PPP^{2})\nabla \epsilon$ be exactly $(\lambda_{s}/\lambda-2\beta)\nabla \epsilon$  to get some cancellation if one follow (the almost exact computation) to handle $R(\epsilon)$ and $(1-\frac{k((\lambda y+\alpha)}{k(\alpha)})|\epsilon|^{2}\epsilon$, we will not go to these part of details, so we are write just like \eqref{eq: ep1s2}).
 \end{remark}
We can now write \eqref{eq: ep1s2} as

\begin{equation}\label{eq: ep1s3}
\begin{aligned}
&i\partial_{s}\epsilon-M(\epsilon)+ib\Lambda \epsilon-i2\beta\nabla \epsilon\\
=&O(\PPP^{2})(y^{2}Q_{\PPP}+yQ_{\PPP}+Q_{\PPP}+i\Lambda Q_{\PPP}+i\nabla Q_{\PPP})\\
+&O(\PPP^{2}i\Lambda \epsilon+O(\PPP^{2})i\nabla \epsilon+R(\epsilon)\\
+&\Psi_{\PPP}+F.
\end{aligned}
\end{equation}

We note that 
\begin{enumerate}
	\item $R(\epsilon)$ and $(1-\frac{k((\lambda y+\alpha)}{k(\alpha)})|\epsilon|^{2}\epsilon$ is not purely perturbative, but can be handled by studying a nonlinear version of energy as illustrated in \cite{raphael2011existence}.
	\item If one views $\Lambda$ as $\sim 1$ (from scaling perspective), and does not care about the extra loss\footnote{In principle, the extra loss of derivative only bothers when one studies $\Re<\beta \nabla \epsilon, \Delta \epsilon>$, which requires that $\epsilon$ in $H^{3/2}$ at first glance, but actually vanishes. } of derivative in $\nabla \epsilon$, then $O(\PPP^{2}i\Lambda \epsilon+O(\PPP^{2})i\nabla \epsilon$  are also of $O(\PPP^{3})$ size, and in particular of lower order compared with $ib\Lambda \epsilon-i2\beta\nabla \epsilon$ on the RHS.
	\item $\Psi_{\PPP}$ terms is large, and cannot be handled perturbatively, this is very different from the case in \cite{raphael2011existence}),it will be treated in the next subsubsection.
\end{enumerate}

Let us neglect $R(\epsilon)$, $O(\PPP^{2}i\Lambda \epsilon+O(\PPP^{2})i\nabla \epsilon$, and $\Psi_{\PPP}$ and revisit the construction of Lyapunov functional in \cite{raphael2011existence} by studying the following (not that simplified) toy model\footnote{We note such a toy model is consistent with the principle mentioned in \cite{raphael2011existence}, which stressing on keep tracking of quadratic terms  in energy type estimates.},
\begin{equation}\label{eq: simplifiedtoy}
\begin{aligned}
&i\partial_{s}\epsilon-M(\epsilon)+ib\Lambda \epsilon-i2\beta\nabla \epsilon\\
=&O(\PPP^{2})(y^{2}Q_{\PPP}+yQ_{\PPP}+Q_{\PPP}+i\Lambda Q_{\PPP}+i\nabla Q_{\PPP})\\
+&F.
\end{aligned}
\end{equation}

Observing $-M+ib\Lambda+i2\beta \nabla$ is indeed self-adjoint with inner product $\Re<f,g>$, it is natural to do energy estimates of \eqref{eq: simplifiedtoy} based on 
\begin{equation}
J(\epsilon):=\Re <M(\epsilon)-ib\Lambda \epsilon+i2\beta \nabla \epsilon,\epsilon>	
\end{equation}
where the Virial term $\int y\partial_{y}\epsilon \bar{\epsilon}$.

Here one crucial technical difficulty will arise since through $\Lambda$ scales like $1$, but the $y\partial_{y}\epsilon$ is not in $H^{1}$, thus one needs to studies a local version of $J$, or $\Lambda$.

Following (3.36) \cite{raphael2011existence}, let $\phi$ be radial, smooth, and monotone in $r$, (or more precisely $|x|$),
\begin{equation}
\begin{cases}		
\phi'(r)=r, r\leq 1\\
\phi'(r)=3-e^{-r}, r\geq 2
\end{cases}
\end{equation}

And let $A$ be large, one define\footnote{$\Lambda_{A}$ is a natural way to truncate $\Lambda$, note that $i\Lambda_{A}$ is still self-adjoint.} $\Lambda_{A}$ as 
\begin{equation}
\Lambda_{A}f:=\frac{1}{2}(\Delta \phi)(\frac{y}{A})f(y)+A(\nabla \phi)(\frac{y}{A})\partial_{y}f
\end{equation}
and let

\begin{equation}
\tilde{J}(\epsilon):=\Re <M(\epsilon)-ib\Lambda_{A} \epsilon+i2\beta \nabla \epsilon,\epsilon>	.
\end{equation}

For certain re-normalization purpose, one will indeed study 
\begin{equation}\label{eq: ly1}
\tilde{I}(\epsilon):=\frac{1}{\lambda^{2}}\Re <M(\epsilon)-ib\Lambda_{A} \epsilon+i2\beta \nabla \epsilon,\epsilon>.	\end{equation}

One has, 
\begin{lemma}\label{lem: lymodel}
Under the bootstrap assumption of Lemma \ref{lem: bootmain}, for \eqref{eq: simplifiedtoy} one has
\begin{equation}\label{eq: lc1}
\lambda^{2}\tilde{I}(\epsilon)\geq \delta_{0}\|\epsilon\|_{H^{1}}^{2}+O(\mathcal \PPP^{3})	
\end{equation}
and 
\begin{equation}\label{eq: mono1}
\frac{d}{dt}\tilde{I}(\epsilon)\equiv \frac{1}{\lambda^{2}}\frac{d}{ds}I_{1}(\epsilon)\geq O(1)
\end{equation}

\end{lemma}
Lemma \ref{lem: lymodel} is essentially contained in \cite{raphael2011existence}, though not in an explicit way.

We will give a quite detailed sketch of proof of Lemma \ref{lem: lymodel} in the last subsubsection of this subsection and we will in particular detailed the discussion of the impact $\beta\nabla \epsilon$ part in the $\tilde{I}$.

And to go back the original full system \eqref{eq: ep1s1}, one will following \cite{raphael2011existence} to study a nonlinear version of $\tilde{I}$ to handle $R(\epsilon)$ and $1-\frac{k(\lambda y+\alpha)}{k(\alpha)}$, i.e. 
\begin{equation}
I(\epsilon):=\frac{1}{\lambda^{2}}\Re<M(\epsilon)-ib\Lambda_{A}\epsilon+i2\beta\nabla \epsilon,\epsilon>-\int \frac{k(\lambda y+\alpha)}{k(\alpha)}F(Q_{\PPP},\epsilon)
\end{equation}
where $F(Q_{\PPP},\epsilon)$ is the real part of all the terms in $\frac{1}{4}|Q_{\PPP}+\epsilon|^{4}$ with order 3 or 4 in $\epsilon$.

And One has 
\begin{equation}\label{eq: lc2}
\lambda^{2}I(\epsilon)\geq \delta_{0}\|\epsilon\|_{H^{1}}^{2}+O(\mathcal \PPP^{3})	
\end{equation}
 and 
 \begin{align}\label{eq: mono2}
 \frac{d}{dt}{I}(\epsilon)& \equiv \frac{1}{\lambda^{2}}\frac{d}{ds}I_{1}(\epsilon)\\ \nonumber  & \geq \frac{1}{\lambda^{4}}\Re <M(\epsilon)-ib\Lambda \epsilon-i2\beta \epsilon,-i\Psi_{\PPP}>+ O(1)\\ \nonumber & =\frac{1}{\lambda^{4}}\Re<M(\epsilon), i\Psi_{\PPP}>+O(1).
 \end{align}
 
 (Since $\Psi_{\PPP}$ is well localised, the difference of $\Lambda_{A}\epsilon$ and $\Lambda \epsilon$ can be neglected when paired with $\Psi_{\PPP}$.)

We note that (up to slight higher order modification), $I(\epsilon)$ is the Lyapounov functional in \cite{raphael2011existence} \textbf{adding an extra term} $\frac{1}{\lambda^{2}}<i2\beta\nabla \epsilon,\epsilon>$. And estimates \eqref{eq: lc2}, \eqref{eq: mono2} has been proved there, (without the $
\beta$ term). \\

The $\beta$ term is of extra smallness in \cite{raphael2011existence} and will not cause a problem, but here we need to put it into the Lyapounov functional, (and it is also natural to do so from the previous discussion).

The handling of $\beta$ term is same as in Lemma \ref{lem: lymodel} and we omit it here. What is more troublesome is the $<M(\epsilon), i\Psi_{\PPP}>$ term. Recall in \cite{raphael2011existence}, one basically has $\Psi_{\PPP}=O(\PPP^{5})$, thus $\frac{1}{\lambda^{2}}<M(\epsilon), \Psi_{\PPP}>$ can be neglected.

However, in the current article, $\Psi_{\PPP}$ is formally $o_{\PPP}(\PPP^{2})$, then a priori $\frac{1}{\lambda^{2}}<M(\epsilon), \Psi_{\PPP}>$ is of size $\frac{1}{\lambda}\sim \frac{1}{t}$, which will cause a logarithm divergence.

We will explain how to handle it by further modifying the Lyapounov functional in the next subsubsection.

\begin{remark}
Such a logarithm divergence	will disappear if one assumes slightly more regularity of $k$, (for example, $C^{2,\alpha}, \alpha>0$), this is also one reason making the $C^{2}$ case subtle.
\end{remark}

\subsubsection{A modified Lyapounov functional}
We define 
\begin{equation}\label{eq: lyapounovi1}
I_{1}(\epsilon):=I+\frac{1}{\lambda^{2}}\Re<-i\epsilon,-i\Psi_{\PPP}>=I+\frac{1}{\lambda^{2}}\Re<\epsilon, \Psi_{\PPP}>.
\end{equation}

And we check $I_{1}$ satisfy Lemma \ref{lem: keyly}

We note that \eqref{eq: lc2} implies 
\begin{equation}
\lambda^{2}I(\epsilon)\geq \delta_{0}\|\epsilon\|_{H^{1}}^{2}+O(\PPP^{3})
\end{equation}
since one can estimate $<\epsilon, \Psi_{\PPP}>$ by $O_{\PPP}(\PPP^{3})$.

It remains to prove
\begin{equation}
\begin{aligned}
&\frac{1}{\lambda^{2}}\partial_{s}(\frac{1}{\lambda^{2}}\Re<-i\epsilon,-i\Psi_{\PPP}>)\\
+&\frac{1}{\lambda^{4}}\Re<M(\epsilon)-ib\Lambda \epsilon-i2\beta\nabla \epsilon, -i\Psi_{\PPP}>\\
\geq &O(1).	
\end{aligned}
\end{equation}

We compute as 
\begin{equation}
\begin{aligned}
&\frac{1}{\lambda^{2}}\partial_{s}(\frac{1}{\lambda^{2}}<-i\epsilon,-i\Psi_{\PPP}>)\\
=&2\frac{b+O(\PPP^{2})}{\lambda^{4}}<\epsilon, \Psi_{\PPP}>\\
+&\frac{1}{\lambda^{4}}<\epsilon,\partial_{s}\Psi_{\PPP}>\\
+&\frac{1}{\lambda^{4}}<-i\epsilon_{s}, i\Psi_{\PPP}>.
\end{aligned}	
\end{equation}
The first two terms on the RHS are both of $O(1)$, (recall also $\partial_{s}\Psi_{\PPP}=o_{\PPP}(\PPP^{3})$, thanks to \eqref{eq: crude}.)

It remains to estimate
\begin{equation}\label{eq: finalcontrol}
<-i\epsilon_{s}, i\Psi_{\PPP}>+<M(\epsilon)-ib\Lambda \epsilon-i2\beta\nabla \epsilon, -i\Psi_{\PPP}>.
\end{equation}
The linear part in the $\epsilon_{s}$ cancels with second part.

Thus \eqref{eq: finalcontrol} are controlled by the $L^{2}$ pair between $\Psi$ and 
\begin{itemize}
\item Terms of size $O(\PPP^{2})Y$, $Y$ is localized and fast decaying.
\item Terms of $O(\PPP^{2})\Lambda \epsilon$, $O(\PPP^{2})\epsilon$
\item Terms of $O(Q\epsilon^{2})$ and $\epsilon^{3}$
\item $\Psi_{\PPP}$ and smaller error $F$	
\end{itemize}
 
All the above pairs are bounded by $O(\PPP^{4})$, thus yielding the desired estimates.

\subsubsection{A sketch of Lemma \ref{lem: lymodel}}
Estimate \eqref{eq: lc1} follows from Lemma \ref{lem: coer}. Note that $M$ is a perturbed version of $L_{1}+iL_{2}$, thus also enjoys Lemma \ref{lem: coer}. The four directions in Lemma \ref{lem: coer} are handled by the orthonality condition, except for $(\epsilon_{1},Q)$, which is bounded by $O(\PPP^{2})$ by mass conservation law, Lemma \ref{lem: conesti}. (The term $<ib\Lambda_{A}\epsilon,\epsilon>, <i2\beta\nabla \epsilon, \epsilon >$ are itself bounded by $O(\PPP^{3})$.)

We now turn to \eqref{eq: lc2}.

One first computes as 
\begin{equation}
\begin{aligned}
\frac{1}{\lambda^{2}}\partial_{s}(\frac{1}{\lambda^{2}}\tilde{I})
=&\frac{b}{\lambda^{4}}\Re<M(\epsilon)-ib\Lambda_{A}\epsilon-i2\beta\epsilon, \epsilon>\\
+&\frac{1}{\lambda^{4}}\Re<M(\epsilon)-ib\Lambda_{A}\epsilon-i2\beta\epsilon, \epsilon_{s}>\\
+&O(1).
\end{aligned}	
\end{equation}
Note that $|b+\frac{\lambda_{s}}{\lambda}|=O(\PPP^{2})$ thanks to Lemma \ref{lem: moduesti}. Also note when the $\partial_{s}$ hits the operator $M_{\PPP}$, $b$, $\beta$ in  $M(\epsilon)-ib\Lambda_{A}\epsilon-i2\beta\epsilon$ those terms becomes lower order and acceptable thanks to \eqref{eq: crude}.

For notation simplicity, let us denote the RHS of \eqref{eq: simplifiedtoy} as $G$. And we note 
\begin{equation}
\Re<M(\epsilon)-ib\Lambda_{A}\epsilon-i2\beta\epsilon, \epsilon_{s}>=<M(\epsilon)-ib\Lambda_{A}\epsilon-i2\beta\epsilon,-i(M\epsilon-ib\Lambda \epsilon+i2\beta \nabla \epsilon)-iG>.	
\end{equation}

It remains to prove
\begin{align}\label{eq: batt}
&\Re<M(\epsilon)-ib\Lambda_{A}\epsilon-i2\beta\epsilon, \epsilon>\\+&\Re<M(\epsilon)-ib\Lambda_{A}\epsilon-i2\beta\epsilon,-i(M\epsilon-ib\Lambda \epsilon+i2\beta \nabla \epsilon)-iG>	\nonumber \\ \geq & O(\PPP^{4}). \nonumber
\end{align}

Note that up to handling the terms involved with $\beta$, estimate \eqref{eq: batt} are  essentially contained in the Lemma 3.3 in \cite{raphael2011existence}. We will sketch the main ideas, highlight the virial algebra, and detailed the computation involving $\beta$.

Let us focus on the second term in \eqref{eq: batt}, which has two parts
\begin{itemize}
\item The part $\Re<M(\epsilon)-ib\Lambda_{A}\epsilon-i2\beta\epsilon,-i(M\epsilon-ib\Lambda \epsilon+i2\beta \nabla \epsilon)>$ which is almost $0$, (and it is indeed zero if one formally taking $A=\infty$, and assuming $\epsilon$ with enough integrability and regularity). 
,but involving a part needs the virial algebra to be estimated combined with the first term of \eqref{eq: batt} with a sign.
\item The part $<M(\epsilon)-ib\Lambda_{A}\epsilon-i2\beta\epsilon,-iG>$ which can be directly estimated by $O(\PPP^{4})$ 
\end{itemize}

For the second part, one first note $-ib\Lambda_{A}\epsilon-i2\beta\epsilon$ are of form $O(\PPP)\epsilon \sim O(\PPP^{2})$, thus a pairing of $G$ gives $O(\PPP^{4})$, (noting the dominate part of $G$ is of form $O(\PPP^{2})$).
For the pairing  $\Re<M(\epsilon), -iG>$, the $F$ part in $G$ is harmless. The pairing $\Re<M(\epsilon), iO(\PPP^{2})((|y|^{2}Q_{\PPP}+yQ_{\PPP}+Q_{\PPP}+i\Lambda Q_{\PPP}+i\nabla Q_{\PPP}))>$ is a prior of size $O(\PPP^{3})$, however, one can use the fact that $M$ be self adjoint, and observe,
\begin{equation}
\begin{aligned}
	&\Re<M(\epsilon), iO(\PPP^{2})(y^{2}Q_{\PPP}+yQ_{\PPP}+Q_{\PPP}+i\Lambda Q_{\PPP}+i\nabla Q_{\PPP})>\\=&O(\PPP^{2})\Re<\epsilon,M(iy^{2}Q_{\PPP}+iyQ_{\PPP}+iQ_{\PPP}i\Lambda Q_{\PPP}-i\nabla Q_{\PPP})>.
\end{aligned}
\end{equation}
Using spectral property \eqref{eq: spectralba1}, \eqref{eq: spectralba2} and plug in orthogonality for $\epsilon$ gives the desired estimate, except for the term $<\epsilon, M(\Lambda yQ)>$ will be estimated by $O(\PPP^{2})+2\Re<\epsilon, Q>$, and the latter is controlled by mass conservation law and estimated also be by $\PPP^{2}$.

For the first part, one first computes 
\begin{align}
&\Re<M(\epsilon)-ib\Lambda_{A}\epsilon-i2\beta\epsilon,-i(M\epsilon-ib\Lambda \epsilon+i2\beta \nabla \epsilon)>	\\
=&\Re <-ib(\Lambda_{A}-\Lambda)\epsilon, -iM(\epsilon)+b\Lambda \epsilon-2\beta \nabla \epsilon>.\nonumber
\end{align}

We need to handle the following three terms $\Re <b(\Lambda_{A}-\Lambda)\epsilon, M(\epsilon)>$, $\Re<ib(\Lambda_{A}-\Lambda)\epsilon,b\Lambda \epsilon>$ and $\Re <b(\Lambda_{A}-\Lambda)\epsilon, -2\beta \nabla \epsilon >.$

As aforementioned, the first term is handled by a crucial Virial algebra in \cite{raphael2011existence}, which we refer for more details. But we mention here one key idea is $<-\Delta \epsilon, \Lambda \epsilon>=2\|\nabla \epsilon\|^{2}$, thus if one view $M$ as $-\Delta$, then (we use the first part in \eqref{eq: batt} $2b\Re <M(\epsilon),\epsilon>+<b(\Lambda_{A}-\Lambda)\epsilon, M(\epsilon)>\approx \int_{|x|\lesssim A} |\nabla \epsilon|^{2}\geq 0$. In practice, we need to use Lemma \ref{lem: coer} to relate $<M(\epsilon), \epsilon>$  and $\|\epsilon\|_{H^{1}}^{2}$, choosing $A$ large to make the truncation of perturbative terms and taking advantage of the exponential decaying of $Q$, so that $M$ can be controlled as $-\Delta+1$ when paired with $b(\Lambda_{A}-\Lambda)\epsilon$.

For the second term, $\Re<ib(\Lambda_{A}-\Lambda)\epsilon,b\Lambda \epsilon>$, one needs to observe since we take the $\Re$, one has 
\begin{equation}
\Re<ib(\Lambda_{A}-\Lambda)\epsilon,b\Lambda \epsilon>=	\Re<ib(\Lambda_{A})\epsilon,b\Lambda \epsilon>=O(b^{2}A^{2}|\epsilon|_{H^{1}}^{2})=O(\PPP^{4}).
\end{equation}

For the last term, one can control $\Re <b(\Lambda_{A}, -2\beta \nabla \epsilon >$ also as $O(b^{2}A|\epsilon|_{H^{1}}^{2})=O(\PPP^{4})$ and one observe the algebra
\begin{equation}
\Re <ib\Lambda \epsilon, \beta \nabla \epsilon>=cb\Re<\beta \nabla \epsilon ,\epsilon>=O(\PPP^{4}).	
\end{equation}

\subsection{Proof of bootstrap estimates}\label{sub: proofofestimate}
\subsubsection{Proof of estimate \eqref{eq: beep}}
Integrate $\partial_{t}I_{1}$ backwards from $t_{0}$ to $t$, $t\in [0,T]$, gives $I_{1}(t)\lesssim t\sim \lambda(t)$, then \eqref{eq: beep} follows from \eqref{eq: lc3}.
\subsubsection{Proof of estimate \eqref{eq: beb}}

Under the assumption \eqref{eq: baep}-\eqref{eq: bab},
estimate \eqref{eq: bab} follows from estimate \eqref{eq: conesti} by conservation laws, when $o_{\PPP(t_{0})}(1)$ is small enough. Note that one has \eqref{eq: atecnicalconstant}, see also Remark \ref{rem: extendconst}.
\subsubsection{Proof of estimate \eqref{eq: belam}}

Taking use of modulation equations \eqref{eq: moduesti} and improved bound \eqref{eq: beb}, we obtain
\begin{align*}
|\frac{d}{dt}(\lambda(t)+\frac{t}{C_{0}})| & =|\lambda_{t}(t)+\frac{1}{C_{0}}|\\
 & =|\lambda_{t}(t)+\frac{b(t)}{\lambda(t)}+\frac{1}{C_{0}}-\frac{b(t)}{\lambda(t)}|\\
 & \leq\delta^{2},
\end{align*}
which implies
\[
|\lambda(t)+\frac{t}{C_{0}}|\leq|\lambda(t_{0})+\frac{t_{0}}{C_{0}}|+\int_{t}^{t_{0}}\delta^{2}ds\leq\delta^{2}\lambda(t).
\]

\subsubsection{Proof of estimate \eqref{eq: bemome}}
WLOG, we consider the first component for
both $\alpha$ and $\beta$, i.e. $\alpha_{1}$ and $\beta_{1}$.
For convenience, we introduce a new time variable $\tau(t)$ with
\[
\tau(t):=\int_{t}^{t_{0}}\frac{1}{\frac{-z}{C_{0}}}dz,
\]

which is equal to
\[
\frac{d\tau(t)}{dt}=-\frac{1}{\frac{-t}{C_{0}}},\;\tau(t_{0})=0.
\]
Then we could transform the modulation equations for $\alpha_{1}$
and $\beta_{1}$ into\footnote{Here we assume $\nabla^{2}k(0)=\begin{pmatrix}-k_{1}\\
 & -k_{2}
\end{pmatrix}$ for some $k_{1},k_{2}>0$, which is a crucial assumption in our analysis during this part.}
\begin{equation}
\frac{d}{d\tau}\begin{pmatrix}\alpha_{1}\\
\beta_{1}
\end{pmatrix}=\begin{pmatrix}0 & -2\\
k_{1} & \frac{1}{C_{0}}
\end{pmatrix}\begin{pmatrix}\alpha_{1}\\
\beta_{1}
\end{pmatrix}+\begin{pmatrix}F_{1}\\
F_{2}
\end{pmatrix},\label{the equation for alpha and beta}
\end{equation}

with the bounds 
\begin{equation}
|F_{1}(t)|+|F_{2}(t)|\lesssim\delta^{2}e^{\frac{\tau(t)}{C_{0}}}\cdotp(-t_{0})\lesssim\delta^{2}\lambda(t),\label{the bound for external force term}
\end{equation}
where we have used the bound for modulation equation \eqref{eq: moduesti} and the improved bounds \eqref{eq: belam}, \eqref{eq: beb}. We have the following technical lemma for $\alpha_{1}$ and $\beta_{1}$,

\begin{lemma}\label{ODE lemma}
Let us consider the ODE system (\ref{the equation for alpha and beta})
with the initial data 
\[
\alpha_{1}(t_{0})=\beta_{1}(t_{0})=0
\]
 and the bound for external force term (\ref{the bound for external force term}),
$\forall t\in[T,t_{0}]$. Then we have the following bound for
$\alpha_{1}$ and $\beta_{1}$,
\[
|\alpha_{1}(t)|+|\beta_{1}(t)|\lesssim\delta^{2}e^{\frac{\tau(t)}{C_{0}}}\cdotp(-t_{0})\lesssim\delta^{2}\lambda(t),\;\forall t\in[T,t_{0}].
\]
\end{lemma}

Applying this lemma in our case, we finish the proof of \eqref{eq: bemome}. Hence the only task is to prove Lemma \ref{ODE lemma}.
\begin{proof}
Step 1. The analysis for the matrix $\begin{pmatrix}0 & -2\\
k_{1} & \frac{1}{C_{0}}
\end{pmatrix}$

Since $k_{1}>0$, we know the eigenvalues of $\begin{pmatrix}0 & -2\\
k_{1} & \frac{1}{C_{0}}
\end{pmatrix}$ are $0<\lambda_{1}\leq\lambda_{2}$, with
\[
\begin{cases}
\lambda_{1}+\lambda_{2}=\frac{1}{C_{0}},\\
\lambda_{1}\lambda_{2}=2k_{1}.
\end{cases}
\]

Up to some rotation, we are faced with the following system,
\begin{equation}
\frac{d}{d\tau}\begin{pmatrix}\tilde{\alpha}_{1}\\
\tilde{\beta}_{1}
\end{pmatrix}=\begin{pmatrix}\lambda_{1}\\
 & \lambda_{2}
\end{pmatrix}\begin{pmatrix}\tilde{\alpha}_{1}\\
\tilde{\beta}_{1}
\end{pmatrix}+\begin{pmatrix}\tilde{F}_{1}\\
\tilde{F}_{2}
\end{pmatrix},\label{the ratation system}
\end{equation}

where $\tilde{F}_{1}$ and $\tilde{F}_{2}$ satisfy the same bound
\eqref{the bound for external force term}. And $\tilde{\alpha}_{1}$, $\tilde{\beta_{1}}$ are non-degenerate
linear combination of $\alpha_{1}$ and $\beta_{1}$. Therefore, we
reduce Lemma \ref{ODE lemma} into proving the same bounds for $\tilde{\alpha}_{1}$
and $\tilde{\beta_{1}}$. 

Step 2. The bounds for $\tilde{\alpha}_{1}$ and $\tilde{\beta_{1}}$

Standard ODE argument tells us, 
\begin{equation}
|\tilde{\alpha}_{1}(\tau)|\leq e^{\lambda_{1}\tau}\int_{0}^{\tau}e^{-\lambda_{1}\tilde{s}}F_{1}(\tilde{s})d\tilde{s},\label{the bound for alpha bar}
\end{equation}

and 
\begin{equation}
|\tilde{\beta}_{1}(\tau)|\leq e^{\lambda_{2}\tau}\int_{0}^{\tau}e^{-\lambda_{2}\tilde{s}}F_{2}(\tilde{s})d\tilde{s}.\label{the bound for beta bar}
\end{equation}

Substituting \eqref{the bound for external force term} into (\ref{the bound for alpha bar}) and (\ref{the bound for beta bar}),
we finish the proof. 
\end{proof}

\begin{remark}
We should highlight that Lemma \ref{ODE lemma} heavily relies on non-degeneracy condition $k_{1}>0$, which means the eigenvalues of $\begin{pmatrix}0 & -2\\
k_{1} & \frac{1}{C_{0}}
\end{pmatrix}$ do not contain zero. In the degenerate case, the matrix $\begin{pmatrix}0 & -2\\
k_{1} & \frac{1}{C_{0}}
\end{pmatrix}$ degenerates in the sense it contains zero eigenvalue. And in this case, $\lambda_{1}=0$ and $\lambda_{2}=\frac{1}{C_{0}}$.

Direct computation only implies,
\begin{align*}
|\tilde{\beta}_{1}(\tau)| & \leq e^{\frac{\tau}{C_{0}}}\int_{0}^{\tau}e^{-\frac{\tau}{C_{0}}}\cdotp\delta^{2}(-t_{0})e^{\frac{\tau}{C_{0}}}d\tau\\
 & \leq e^{\frac{\tau}{C_{0}}}\cdotp\delta^{2}(-t_{0})\tau\\
 & \leq\delta^{2}(-t)\ln\frac{t}{t_{0}},
\end{align*}

which is not enough to close the estimate \eqref{eq: bemome}.

\end{remark}

\begin{remark}
Under radial assumption, $\alpha(t)=\beta(t)=0$, and \eqref{eq: bemome} automatically holds. Therefore, we could drop the assumption $det\nabla^{2}k(0)\neq 0$ in \eqref{as: main} for this special case.
\end{remark}
\subsection{Proof of Lemma \ref{lem: conesti}, \ref{lem: moduesti}}\label{subsec: tech}
The proof of Lemma \ref{lem: conesti}, \ref{lem: moduesti} are similar to the proof of Lemma 3.1 in \cite{raphael2011existence}. We do a quite detailed sketch for the convenience of the readers.

\subsubsection{Proof of Lemma \ref{lem: conesti}}
We first note that here the key algebra is 
\begin{equation}\label{eq: alge1}
E(Q+\epsilon)+\frac{1}{2}M(Q+\epsilon)=E(Q)+M(Q)+O(\epsilon^{2}).	
\end{equation}
The first order term in $\epsilon$ vanishes due to \eqref{eq: groundstate}.

We now go to more more details for \eqref{eq: conesti}.

One notes by\eqref{eq: massini} and 
\begin{equation}
M(u)=M(u_{0})=\|Q\|_{2}^{2}+O(\lambda^{4}),	
\end{equation}
Meanwhile, 
\begin{equation}
k(\alpha)M(u)=|Q_{p}+\epsilon|_{L_{2}}^{2}.	
\end{equation}
Plug in $k(\alpha)=1+\frac{1}{2}\nabla^{2}k(0)(\alpha,\alpha)$ and \eqref{eq: qbmass}, one obtains
\begin{equation}\label{eq: massesti}
-\frac{\nabla^{2}k(0)(\alpha,\alpha)}{2}Q^{2}+2\Re<Q_{\PPP},\epsilon>+O(\|\epsilon\|_{H_{1}}^{2})=o_{\PPP}(P^{2}).	
\end{equation}

Also, one has 
\begin{equation}
\begin{aligned}
&k(\alpha)\lambda^{2}\tilde{E}(u_{0})=\frac{1}{2}\int |\nabla Q_{\PPP}+\nabla \epsilon|^{2}-\frac{1}{4}\int \frac{k(\lambda y+\alpha)}{k(\alpha)}|Q_{\PPP}+\epsilon|^{4}\\
=&\tilde{E}(Q_{\PPP})+\Re<\epsilon, -\Delta Q_{\PPP}-\frac{k(\lambda y+\alpha)}{k(\alpha)}|Q_{\PPP}|^{2}Q_{\PPP}>+O(\|\epsilon\|_{H^{1}}).	
\end{aligned}
\end{equation}
Plug in \eqref{eq: qbenergy}, (and do Taylor expansion for $k$ up to second order), one obtains
\begin{equation}\label{eq: energyesti}
\begin{aligned}
\lambda^{2}E_{in}(u_{0})=&\frac{b^{2}}{8}\int|y|^{2}Q^{2}+\frac{|\beta|^{2}}{2}\int Q^{2}-\frac{\lambda^{2}}{8}\int\nabla^{2}k(0)(y,y)Q^{4}\\
+&\Re<\epsilon, -\Delta Q_{\PPP}-\frac{k(\lambda y+\alpha)}{k(\alpha)}|Q_{\PPP}|^{2}Q_{\PPP}>
  +O(\|\epsilon\|_{H^{1}}^{2})+o_{\PPP}(\PPP^{2}).	
\end{aligned}
\end{equation}

Combing \eqref{eq: massesti} and \eqref{eq: energyesti}, one obtains
\begin{equation}
\begin{aligned}
  &\frac{b^{2}}{8}\int|y|^{2}Q^{2}+\frac{|\beta|^{2}}{2}\int Q^{2}-\frac{\nabla^{2}k(0)(\alpha,\alpha)}{4}\int Q^{2}-\lambda^{2}(E_{0}+\frac{1}{8}\int\nabla^{2}k(0)(y,y)Q^{4})\\
  =&-\Re <\epsilon, -\Delta Q_{\PPP}+Q_{\PPP}\frac{k(\lambda y+\alpha)}{k(\alpha)}|Q_{\PPP}|^{2}Q_{\PPP}>+O(\|\epsilon\|_{H^{1}}^{2})+o_{\PPP}(\PPP^{2}).
  \end{aligned}
\end{equation}

It remains to bound $\Re <\epsilon, -\Delta Q_{\PPP}+Q_{\PPP}\frac{k(\lambda y+\alpha)}{k(\alpha)}|Q_{\PPP}|^{2}Q_{\PPP}>$, we apply \eqref{eq: qb}, and \textbf{crucially} use orthogonality condition to cancel the $O(\PPP)$ terms $ib\Lambda Q_{\PPP}, -2i\beta \nabla Q_{\PPP}$ in \eqref{eq: qb}, we obtain
\begin{equation}
	\Re <\epsilon, -\Delta Q_{\PPP}+Q_{\PPP}\frac{k(\lambda y+\alpha)}{k(\alpha)}|Q_{\PPP}|^{2}Q_{\PPP}>=O(\PPP^{2}\|\epsilon\|_{H^{1}}).
\end{equation}

Lemma \ref{lem: conesti} thus follows.

\subsubsection{Proof of Lemma \ref{lem: moduesti}}

 This part is totally parallel with the content in Appendix A for \cite{raphael2011existence}, and the only difference is the size of error caused by $C^{2}$ property of $k(x)$. For reader's convenience, let us state it as follows. 

Step 1. The estimate for the linear part of $\epsilon(t)$ equation \eqref{eq: ep1s3}

This step is exactly the same with (A.3)-(A.7) of Appendix A in \cite{raphael2011existence}. We compute the inner products between the linear part of $\epsilon(t)$ equation \eqref{eq: ep1s3} and some special directions, and obtain
\begin{equation}\label{linear for gradQ}
\Re<-M(\epsilon)+ib\Lambda \epsilon-i2\beta\nabla \epsilon, \nabla Q_{\PPP}>=O(\PPP^{2}\|\epsilon\|_{H^{1}}),
\end{equation}

\begin{equation}\label{linear for yQ}
\Im<-M(\epsilon)+ib\Lambda \epsilon-i2\beta\nabla \epsilon, y Q_{\PPP}>=2\beta (\epsilon_{1},Q)+O(\PPP^{2}\|\epsilon\|_{H^{1}}),
\end{equation}

\begin{equation}\label{linear for lambdaQ}
\Re<-M(\epsilon)+ib\Lambda \epsilon-i2\beta\nabla \epsilon, \Lambda Q_{\PPP}>=-2\Re<\epsilon,Q_{\PPP})+O(\PPP^{2}\|\epsilon\|_{H^{1}}),
\end{equation}

\begin{equation}\label{linear for y2Q}
\Im<-M(\epsilon)+ib\Lambda \epsilon-i2\beta\nabla \epsilon, |y|^{2} Q_{\PPP}>=O(\PPP^{2}\|\epsilon\|_{H^{1}}),
\end{equation}
\begin{equation}\label{linear for rho}
\Re<-M(\epsilon)+ib\Lambda \epsilon-i2\beta\nabla \epsilon, \rho e^{-i\frac{b|y|^{2}}{4}+i\beta y}> =O(\PPP^{2}\|\epsilon\|_{H^{1}}).
\end{equation} 
  
Step 2. The derivation for the modulation equations

Now we take the full $\epsilon(t)$ equation \eqref{eq: ep1s3} with the above five directions, and using mass conservation law derive,
\begin{align}\label{the law for b}
 & -\frac{||yQ||_{2}^{2}}{4}(b_{s}+b^{2}-d_{0}(\alpha,\alpha)) \\
= & o_{\PPP}(\PPP^{2})+|Mod(t)|(||\varepsilon||_{H^{1}}+\PPP)+||\varepsilon||_{H^{1}}^{2}+\PPP^{2}||\varepsilon||_{H^{1}},\nonumber 
\end{align}
\begin{align}\label{the law for lamuda}
 & -(\frac{\lambda_{s}}{\lambda}+b)\int\Lambda Q\cdotp|y|^{2}Q \\
= & o_{\PPP}(\PPP^{2})+|Mod(t)|(||\varepsilon||_{H^{1}}+\PPP)+||\varepsilon||_{H^{1}}^{2}+\PPP^{2}||\varepsilon||_{H^{1}},\nonumber 
\end{align}
\begin{align}\label{the law for alpha}
 & -\int[(\frac{\alpha_{s}}{\lambda}-2\beta)\cdotp\nabla Q]\cdotp yQ \\
= & o_{\PPP}(\PPP^{2})+|Mod(t)|(||\varepsilon||_{H^{1}}+\PPP)+||\varepsilon||_{H^{1}}^{2}+\PPP^{2}||\varepsilon||_{H^{1}},\nonumber 
\end{align}
\begin{align}\label{the law for beta}
 & -\int[(\beta_{s}+b\beta-c_{0}(\alpha)\lambda)\cdotp yQ]\cdotp\nabla Q\\
= & o_{\PPP}(\PPP^{2})+|Mod(t)|(||\varepsilon||_{H^{1}}+\PPP)+||\varepsilon||_{H^{1}}^{2}+\PPP^{2}||\varepsilon||_{H^{1}},\nonumber 
\end{align}
\begin{align}\label{the law for the phase}
 & (\tilde{\gamma}_{s}-|\beta|^{2}+d_{1}(\alpha,\alpha))\int Q\rho \\
= & o_{\PPP}(\PPP^{2})+|Mod(t)|(||\varepsilon||_{H^{1}}+\PPP)+||\varepsilon||_{H^{1}}^{2}+\PPP^{2}||\varepsilon||_{H^{1}}.\nonumber 
\end{align}
Under the assumption $||\varepsilon||_{H^{1}}+|\PPP|<<1$ and \eqref{the law for b}-\eqref{the law for the phase}, we finish the proof of Lemma \ref{lem: moduesti}.  

\begin{remark}
In this part, our analysis only changes an error of size $O(\PPP^{3})$ compared with \cite{raphael2011existence}. Therefore, the above estimates could be derived by the Appendix A in \cite{raphael2011existence}.

\end{remark}
  

 \bibliographystyle{amsplain}
\bibliography{BG.bib}
\bibliographystyle{plain}
\end{document}